\documentclass[12pt]{article}
\usepackage{latexsym}
\usepackage{amssymb}
\usepackage{epsfig,graphics,psfrag}
\usepackage{geometry}                % See geometry.pdf to learn the
\usepackage{graphicx}
\usepackage{latexsym}
\usepackage{amsmath}
\usepackage{float}

\def\eqn#1{\def\theequation{#1}}

\def\doublespaced{\baselineskip=\normalbaselineskip
    \multiply\baselineskip by 2}
\def\doublespace{\doublespaced}
\def\singlespaced{\baselineskip=\normalbaselineskip}
\def\singlespace{\singlespaced}
\newcommand{\hd}{\mbox{H-dim}\;}
\doublespace

\newcommand{\rn}[1]{\mathbb{R}^{#1}}

\newcommand{\re}{ \mathbb{R}}
\newcommand{\ce}{ \mathbb{C}}
\newcommand{\ima}{\mbox{ Im }} 
\newcommand{\rea}{\mbox{ Re }} 
\newcommand{\Hh}{ \mathbb{H}}
\newcommand{\beq}{\begin{equation}}
\newcommand{\bea}[1]{\begin{array}{#1} }
\newcommand{\eeq}{ \end{equation}}
\newcommand{\ea}{ \end{array}}
\newcommand{\ep}{\epsilon}
\newcommand{\es}{\emptyset}

\newcommand{\al}{\alpha}
\newcommand{\defeq}{ := }
\newcommand{\ga}{\gamma}

\newcommand{\de}{\delta}
\newcommand{\ds}{\displaystyle}
\newcommand{\ts}{\textstyle}
\newcommand{\rar}{\mbox{$\rightarrow$}}

\newcommand{\ran}{\rangle}
\newcommand{\lan}{\langle}
\newcommand{\Ga}{\Gamma}
\newcommand{\la}{\lambda}
\newcommand{\La}{\Lambda}
\newcommand{\ar}{\partial}
\newcommand{\si}{\sigma}

\newcommand{\om}{\omega}
\newcommand{\Om}{\Omega}
\newcommand{\be}{\beta}

\newcommand{\ph}{\phi}
\newcommand{\he}{\theta}
\newcommand{\He}{\Theta}

\newcommand{\hs}[1]{\mbox{$ \hspace{#1}$}}

\newcommand{\sem}{\setminus}
\newcommand{\ze}{\zeta}

\newcommand{\ti}{\tilde}

\def\singlespace{\singlespaced}
\singlespace  

\begin{document}

\title{$p$ Harmonic Measure in Simply Connected Domains}
\author{John L. Lewis\thanks{%
email: john@ms.uky.edu} \thanks{%
Lewis was partially supported
by DMS-0552281.}\linebreak \\
Department of Mathematics, University of Kentucky\linebreak \\  Lexington, KY  40506-0027, USA\\
Kaj Nystr\"{o}m\thanks{%
email: kaj.nystrom@math.umu.se} \thanks{%
Nystr\"{o}m was  partially supported by grant VR-70629701 from the Swedish research council VR.}\linebreak \\
 Department of Mathematics, Ume{\aa} University\linebreak \\ S-90187 Ume{\aa}, Sweden\\
 Pietro Poggi-Corradini\thanks{%
email: pietro@math.ksu.edu}\linebreak \\
 Department of Mathematics, Cardwell Hall, Kansas State University \linebreak \\
Manhattan, KS 66506, USA}
\maketitle

\begin{abstract}
\noindent Let $ \Om $ be a  bounded simply connected domain in the complex plane, $ \ce  $. Let  $  N $ be a 
neighborhood of $ \ar \Om$, let  $ p $ be fixed,  $ 1 < p < \infty, $ and let
 $ \hat  u $ be a positive weak solution   to the $ p $ Laplace equation
in $  \Om \cap N.  $  Assume that  $ \hat  u $ has  zero  boundary values on  
$ \ar  \Om$ in the  Sobolev sense and extend $ \hat  u $ to  $  N \sem  \Om  $ by putting 
$ \hat  u \equiv 0 $ on $ N  \sem \Om. $   Then  there exists a positive finite Borel measure $ \hat  \mu $  on 
$ \ce  $ with support contained in $ \ar \Om $ and such
that 
\begin{eqnarray*}   
 \int | \nabla \hat  u |^{p - 2}  \, \lan  \nabla \hat  u ,   
 \nabla \phi  \ran \, dA  = -  \int  \phi  \, d \hat  \mu  
 \end{eqnarray*} 
 whenever  
 $ \ph \in C_0^\infty (  N ). $  If $ p = 2$ and if $\hat u$ is 
the Green function for $\Omega$ with pole at $x\in\Omega\setminus\bar N$ then the measure $\hat\mu$ coincides with harmonic measure at 
$x$, $\omega=\omega^x$, associated to the Laplace equation. In this paper we continue the studies 
in [BL05], [L06] by establishing new results, in simply connected domains, concerning the Hausdorff dimension of the support
of the measure $\hat\mu$. In particular, we prove results, for $ 1 < p < \infty$, $p\neq 2$, reminiscent of the famous 
result of Makarov [Mak85] concerning the Hausdorff dimension of the support of harmonic measure in simply connected domains.\\


\end{abstract}

\newpage

\section{Introduction}   

Let $\Om\subset\mathbf R^{n}$ be a bounded domain and recall that the
continuous Dirichlet problem for Laplace's equation in $\Omega$ can be
stated as follows. Given a continuous function $f$ on $\partial\Omega$,
find a harmonic function $u$ in $\Omega$ which is continuous in
$\overline{\Omega}$, with $u=f$ on $\partial\Omega$. Although such a
classical solution may not exist, it follows from a method of
Perron-Wiener-Brelot that there is a unique bounded harmonic function
$u$ with continuous boundary values equal to $f$, outside a
set of capacity zero (logarithmic capacity for $n=2$ and Newtonian
capacity for $n>2$). The maximum principle and Riesz 
representation theorem yield, for each $x\in\Omega$, the existence of
a Borel measure $\omega^{x}$ with $\omega^{x} (\partial\Om)=1,$ and 
\begin{equation*}
u (x)=\int_{\partial\Omega}f (y)d\omega^{x} (y)\quad \mbox{whenever $x\in \Om$.}
\end{equation*}
Then, $\om =\om^{x}$ is referred to as the harmonic measure at $x$
associated with the Laplace operator.

Let 
also $g=g(\cdot)=g(\cdot,x)$ be the Green function for $\Omega$ with pole at $x\in\Omega$ and extend 
$g$ to $\mathbf R^n\setminus\Omega$ by putting $g\equiv 0$ on $\mathbf R^n\setminus\Omega$. Then $\omega$ is the Riesz measure associated to $g$ in the sense that 
\begin{eqnarray*}\label{intro1}
{ \ds \int }  \lan \nabla g, \nabla \ph \ran \, dx=  -  { \ds \int }   \, \ph \, d \omega\mbox{ whenever  
$\ph \in C_0^\infty (  \mathbf R^n\setminus\{x\})$}.
   \end{eqnarray*}
We define the  
Hausdorff dimension of $ \omega  $, denoted    $\hd  \omega$,  by 
\[    \hd  \omega \, =
 { \ds \inf }  \{ \alpha : \mbox{ there exists $E$ Borel $ 
 \subset \ar \Om $ with }
  H^\alpha ( E ) = 0 \mbox{ and }    \omega ( E ) =  \omega ( \ar \Om ) \}, 
  \]
  where $H^\alpha(E)$, for $\alpha\in\mathbf R_+$, is the
$\alpha$-dimensional Hausdorff measure of $E$ defined below. In the past twenty years a number of 
  remarkable results concerning  $\hd  \omega$ have been established in planar domains, 
  $\Omega\subset\mathbf R^2$. In particular, Carleson [C85] showed that $\hd \omega = 1$
when $\partial\Omega$ is a snowflake and that $\hd \omega \leq 1$
for any  self similar  Cantor set. Later
Makarov [Mak85] proved that $\hd \omega = 1$ for any simply connected domain in the plane. Furthermore, Jones 
and Wolff
[JW88] proved that $\hd \omega \leq 1$ whenever $\Omega\subset\mathbf R^2$ and $\omega$ exists and 
Wolff [W93] strengthened [JW88] by showing
that $\omega$ is concentrated on a set of s finite $H^1$-measure. We also mention results
of Batakis [Ba96], Kaufmann-Wu [KW85], and Volberg [V93] who showed, for certain fractal domains and
domains whose complements are Cantor sets, that
\begin{eqnarray*}
\mbox{Hausdorff dimension of $\partial\Omega$ = $\inf\{\alpha:\  H^\alpha(\partial\Omega)=0\}>
\hd \omega$.}
\end{eqnarray*}
Finally we note that higher dimensional results for the dimension of harmonic measure can be found
in [Bo87], [W95], and [LVV05].

In [BL05] the first author, together with Bennewitz, started the study of the dimension of a measure, here 
referred to as $p$ harmonic measure, associated with a 
positive $p$
harmonic function which vanishes on the boundary of certain domains
in the plane. The study in [BL05] was continued in  
[L06]. Let $ \ce $ denote the complex plane and let $ dA  $ be  Lebesgue measure on    $ \ce.  $
If $ O\subset\ce $ is open and 
$ 1  \leq  q  \leq  \infty, $ let  $ 
W^{1,q} ( O ) $ be the space of equivalence classes of functions
$ \hat   u $ with distributional gradient $  \nabla \hat   u = ( \hat   u_{x}, 
 \hat  u_{y} ), $ both of which are $ q $ 
th power integrable on $ O. $  Let 
\[ \| \hat  u \|_{1,q} = \| \hat  u \|_q +  \| \nabla \hat  u \|_{q}   \]
be the  norm in $ W^{1,q} ( O ) $ where $ \| \cdot \|_q $ denotes
the usual  Lebesgue $ q $ norm in $ O. $  Let $ C^\infty_0 (O )$ be
 infinitely differentiable functions with compact support in $
O $ and let  $ W^{1,q}_0 ( O ) $ be the closure of $ C^\infty_0 ( O ) $ 
in the norm of $ W^{1,q} ( O  ). $   
Let $ \Om \subset \ce $ be a  simply connected domain and suppose that 
the  boundary of $  \Om$, $  \ar \Om$,  
is bounded and non empty.  Let  $  N $ be a 
neighborhood of $ \ar \Om, \,   p $ fixed,  $ 1 < p < \infty, $ and let
 $ \hat  u $ be a positive weak solution   to the $ p $ Laplace equation
in $  \Om \cap N.  $ That is, 
$ \hat  u \in W^ {1,p} ( \Om \cap N ) $ and  
\eqn{1.1} \beq \int | \nabla \hat  u |^{p - 2}  \, \lan  \nabla \hat  u ,   
 \nabla \he  \ran \, dA   = 0    \eeq 
 whenever $  \he  \in 
W^{1, p}_0 (  \Om\cap N). $    
 Observe that if $ \hat  u $ is smooth and $ \nabla \hat  u \not = 0 $ 
in $ \Om \cap N, $ then   
$ \, \nabla \cdot ( | \nabla \hat  u |^{ p - 2} \, \nabla \hat  u ) \equiv 0,
$  in the classical sense, where $ \nabla \cdot $ denotes divergence.  We assume that  $ \hat  u $ has  zero  boundary values on  
$ \ar  \Om$ in the  Sobolev sense. More specifically if 
$ \ze \in C^\infty_0 (   N ), $ then 
$ \hat  u \, \ze \in W^{1,p}_0 (  \Om\cap N). $ 
Extend $ \hat  u $ to  $  N \sem  \Om  $ by putting 
$ \hat  u \equiv 0 $ on $ N  \sem \Om. $   Then  $ \hat  u \in W^{1,p}
(  N )  $   and it follows from (1.1), as in [HKM93], that 
there exists a positive finite Borel measure $ \hat  \mu $  on 
$ \ce  $ with support contained in $ \ar \Om $ and the property
that 
\eqn{1.2} \beq    
 \int | \nabla \hat  u |^{p - 2}  \, \lan  \nabla \hat  u ,   
 \nabla \phi  \ran \, dA  = -  \int  \phi  \, d \hat  \mu   \eeq 
 whenever  
 $ \ph \in C_0^\infty (  N ). $ 
 We note that if $ \ar \Om $ is smooth enough, then 
$ \,  \,  d \hat  \mu = | \nabla \hat  u |^{p - 1} \, d H^1 |_{\ar \Om}. $ Note that if $ p = 2$ and if $\hat u$ is 
the Green function for $\Omega$ with pole at $x\in\Omega$ then the measure $\hat\mu$ coincides with harmonic measure at $x$, $\omega=\omega^x$, 
introduced above. We refer to $\hat\mu$ as the $p$ harmonic measure associated to $\hat u$. In [BL05], [L06] the 
Hausdorff dimension of the $p$ harmonic measure $\hat\mu$ is studied for general $p$, $1<p<\infty$, and 
to state results from [BL05], [L06]  we next properly introduce the notions of Hausdorff measure and Hausdorff dimension. 
In particular, let 
points in the 
complex plan  be denoted by $ z = x + i y $ and put 
 $ B (z, r ) = 
\{ w   \in \ce  : |  w  -  z | < r \} $ whenever $  z  \in \ce $
and $ r > 0. $    Let $ d ( E, F ) $  denote the distance between
the sets $ E, F \subset \ce $. If $ \la > 0 $ is a positive function on $ (0,r_0 ) $   with 
 $ { \ds \lim_{r\rar 0} \la (  r ) = 0 } $ define  $ H^\la $
Hausdorff measure on   $ \ce $       
as follows:  For fixed $0 <  \delta
 < r_0  $ and $ E \subseteq \rn{2} $, let
$ L ( \delta ) = \{ B ( z_i,  r_i ) \} $ be such that
$ E \subseteq \bigcup \, B ( z_i , r_i ) $ and $ 0 < r_i < \delta , ~~
i = 1,2,..$.  Set
\[ \phi_{ \delta }^\la (E) = {\displaystyle \inf_{ L ( \delta ) } }
 \sum  \,  \la ( r_i )  .  \]
  Then
\[ H^\la (E) = {\displaystyle \lim_{ \delta \rightarrow 0 } }
 \, \, \phi_{ \delta }^\la (E)  . \]
In case $ \la ( r ) = r^\alpha $ we write $ H^\alpha $ for $  H^\la. $ We now define the  
Hausdorff dimension of the measure $  \hat \mu  $ introduced in (1.2) as 
\[    \hd  \hat \mu  \, =
 { \ds \inf }  \{ \alpha : \mbox{ there exists $E$ Borel $ 
 \subset \ar \Om $ with }
  H^\alpha ( E ) = 0 \mbox{ and }    \hat \mu ( E ) =  \hat \mu ( \ar \Om ) \}.  
  \]
  
In [BL05] the first author, together with Bennewitz, proved the following theorem. \\ 

\noindent {\bf Theorem A.} { \em Let $ \hat  u, \hat  \mu,  $ be as in (1.1),
(1.2). If $   \ar  \Om   $ is   a quasicircle,  then 
 $  \hd \hat  \mu \leq  1 $ for  $  2  \leq  p < \infty ,  $  while  
$ \hd \hat  \mu \geq  1$ 
for $   1 < p \leq  2. $ Moreover, if  $ \ar \Om $ is the von
Koch snowflake then strict inequality holds for  $ \hd \hat  \mu. $ } \\

In [L06] the  results in [BL05] were improved at the expense of assuming
more about $ \ar \Om$. In particular, we refer to [L06] for the definition of a  $ k $ quasi-circle. The following 
theorem is proved in [L06]. \\

\noindent {\bf Theorem B.}
 { \em Given  $ p,  1 < p  < \infty, p \not
= 2, \,  $
 there exists $  k_0 ( p ) > 0$  such that  if  $ \ar \Om $
 is a $ k $ quasi-circle and  $ 0 < k < k_0 ( p), $
 then
 \[  \bea{l} (a) \hs{.2in}    \hat  \mu  \mbox{ is concentrated on a set of $ \si $
finite $ H^1 $ measure when }    p > 2.
\\
(b) \hs{.2in} \mbox{There exists $  A =   A ( p ),  0 <  A( p  ) <
\infty, $ such
that if $ 1   < p < 2,  $ then  $ \hat  \mu $ is absolutely } \\
 \hs{.38in} \mbox{  continuous
with respect  to  Hausdorff measure defined  relative  to  $
  \ti \la $ where }   \ea \] }    \[   \ti   \la ( r ) =
r \, \exp[   A  \sqrt{ \log 1/r \,  \log \log \log 1/r } ],
0 < r < 10^{- 6}.   \]\

We note  that  Makarov in [Mak85] proved Theorem B  for
harmonic measure $\omega$, $ p = 2 $,  when  $
    \Om    $ is simply   connected.  Moreover, in this case it suffices to take 
 $  A  =  6 \sqrt{ (\sqrt{24} - 3)/{5}}$, see [HK07].  
   In this paper we continue the studies in [BL05] and [L06] and we prove the following theorem.\\

\noindent {\bf Theorem 1.}
 { \em Given  $ p,  1 < p  < \infty, p \not
= 2, \,  $ let  $ \hat  u, \hat  \mu  $ be as in (1.1), (1.2), 
 and  suppose 
$ \Om $ is simply connected. 
Put  \[   \la ( r ) =
r \, \exp[   A  \sqrt{ \log 1/r \,   \log \log 1/r } ],
0 < r < 10^{- 6}.   \] Then the following is true.  
 \[  \bea{l} (a) \hs{.2in}     \mbox{  If   $ p > 2,  $    there exists 
$ A = A (
p ) \leq  - 1 $ such that $ \hat  \mu $ is concentrated } 
\\ \hs{.38in}  \mbox{  on a set of $ \si $
finite $ H^\la  $ measure.  }    
\\
(b) \hs{.2in} \mbox{  If  $   1 < p < 2, $  
there exists $  A =   A ( p ) \geq 1, $ 
 such
that  $  \hat \mu $ is absolutely } \\
 \hs{.38in} \mbox{  continuous
with respect  to   $ H^\la. $ }   \ea \] }       

Note  that Theorem 1  and the definition of  $ \hd \hat \mu $
imply  the following corollary. \\

\noindent {\bf Corollary 1.} {\em Given  $ p,  1 < p  < \infty, p \not
= 2, \,  $ let  $ \hat  u, \hat  \mu  $ be as in (1.1), (1.2), 
 and  suppose 
$ \Om $ is simply connected. Then 
 $  \hd \hat  \mu \leq  1 $ for  $  2  \leq  p < \infty ,  $  while  
$ \hd \hat  \mu \geq  1$ 
for $   1 < p \leq  2. $}  
\\

In  Lemma 2.4, stated below,  we first show
that it is enough to to prove Theorem 1 for a specific $p$  harmonic
function $ \hat  u $ satisfying the hypotheses. Thus, we choose $ z_0 \in
\Omega  $ and let $ u $ be the $ p $ capacitary functions for 
$   D  = \Om \sem \overline{B} ( z_0, d ( z_0, \ar \Om )/2). $  Then $ u $
is  $ p $ harmonic in $ D $ with continuous boundary values, $ u
\equiv 0 $ on $ \ar  \Om  $ and $ u \equiv 1 $ on $ \ar B ( z_0, 
d ( z_0, \ar \Om )/2 ).$ Furthermore, to prove Theorem 1, we build on the tools 
and techniques developed in [BL05]. In particular, as noted in 
[BL05, sec.~7, Closing Remarks, problem 5], given the tools in [BL05] the main difficulty in establishing Theorem 1 is to prove the following 
result.\\

\noindent {\bf Theorem 2.} {\em Given  $ p,  1 < p  < \infty, p \not
= 2, \,  $ let $
 u, D $  be as above.  There exists  $ c_1 \geq 1, $ depending
only on $  p$,  such that  
\[
c_1^{-1}\frac{u (z)}{d (z,\ar \Om)}\leq |\nabla u (z)|\leq
c_1 \frac{u (z)}{d (z,\ar \Om)},  \mbox{ whenever $ z \in D.$ } 
\] }

In fact, most of our effort in this paper is devoted to proving
Theorem 2.  Armed with Theorem 2 we then use  arguments from
[BL05]  and additional measure-theoretic arguments to prove
Theorem 1.  To further appreciate and understand the importance of the type of estimate we establish in Theorem 2, we note that 
this type of estimate is also crucial in the recent work by the first and second author on the boundary behaviour, regularity and 
free boundary regularity for 
$p$ harmonic functions, $p\neq 2$, $1 <   p < \infty$, in domains in $\mathbf R^n$, $n\geq 2$, which 
are Lipschitz or which 
are well approximated by Lipschitz domains in 
the Hausdorff distance sense, see [LN07,LN,LN08a,LN08b]. Moreover, Theorem 2 seems likely 
to be an important step when trying to solve  several problems for $p$ harmonic functions and $p$ harmonic measure, in planar 
simply-connected domains previously only studied in the case $p=2$, i.e., for harmonic functions and harmonic measure. 
In particular, we refer to [BL05, sec.~7, Closing Remarks] and [L06, Closing Remarks] for discussions of open problems.

The rest of the paper is organized as follows. In section 2 we list
some basic local  results for a positive $ p $ harmonic function vanishing
on a portion of $ \ar \Om. $  In section 3 we use these results
to prove Theorem 1 under the assumption that  Theorem 2 is valid. 
In sections  4 and 5  we then prove Theorem 2.  
  
Finally the first author would like to thank Michel Zinsmeister
for some helpful comments regarding the proof of (4.16).

\noindent   \section{Basic Estimates.} 
In  the sequel  $ c $ will denote a  positive constant  $ \geq 1  $ (not
necessarily the same at each occurrence), which may depend only on $ p, $ unless otherwise stated.
In general, $ c ( a_1, \dots, a_n ) $ denotes a positive constant 
$ \geq 1, $  which may depend only on $ p, a_1, \dots, a_n, $ not
necessarily the same at each occurrence.  $ C $ will denote an
absolute constant.  $ A \approx B $ means
that  $  A/B $ is bounded above and below by positive constants
depending only on $ p. $ 
 In this section, we  will   always assume  that $ \Om $ is a
 bounded  simply connected domain,  $ 0 < r  < \mbox{ diam } \ar
\Om$   and $ w \in \ar \Om$.  We begin by stating  some 
interior and boundary estimates for $ \ti u, $  a positive weak solution to 
the $ p $ Laplacian  in $ B ( w,  4 r ) \cap \Om $  with  
$ \ti u  \equiv 0 $ in the Sobolev sense 
  on $ \ar \Om \cap B ( w, 4 r ).$ 
 That is, 
 $ \ti u \in W^{ 1, p} ( B ( w, 4 r ) \cap \Om ) $ and (1.1) holds whenever 
 $ \he \in  W^{1, p}_0 ( B ( w, 4 r ) \cap \Om ). $  Also  $ \ze \ti u \in 
 W^{ 1, p }_0  (  B ( w, 4 r )  \cap \Om   ) $ whenever $ \ze \in 
C_0^\infty ( B ( w, 4 r ) ). $  Extend $ \ti u $ to $ B ( w, 4 r ) $
by putting $ \ti u \equiv 0 $ on $ B ( w, 4 r )\sem \Om. $ Then there exists 
a locally finite positive Borel measure $ \ti \mu $ with 
support $ \subset B ( w, 4 r ) \cap \ar \Om $  and  for which  
(1.2) holds with $ \hat u $ replaced by $ \ti u $  and $ \ph 
\in C_0^\infty ( B ( w, 4 r ) ). $ 
  Let $ { \ds \max_{B ( z, s )}  \ti u, \,  \min_{B ( z, s)} \ti u } $ be the 
 essential supremum and infimum  of $ \ti u $ on $ B ( z, s) $ 
whenever $ B ( z, s ) \subset B ( w, 4 r ). $  
 For references to proofs of Lemmas 2.1 - 2.3 (see [BL05]). \\ 

\noindent  {\bf Lemma 2.1.} {\em  Fix $ p, 1 < p < \infty, $ and
let $\Om,  w, r, \ti u,   $ be as above.  
   Then  \[ c^{ - 1}  r^{ p - 2} \, \int_{B ( w, r/2)}  \,  | \nabla \ti u 
|^{ p } \,  dx  \, \leq \, \max_{ B ( w, r ) } \, \ti u^p  \, \, \leq 
c \, r^{ - 2} \, \int_{B ( w, 2 r ) } \, \ti u^p \, dx.  \]  If 
$ B ( z,   2s  ) \subset \Om, $ then } 
\[ \max_{B ( z , s  ) } \, \ti u  \, \leq  c \min_{ B ( z , s  )} \ti u. \] 
\\

\noindent 
{\bf Lemma 2.2. } {\em  Let $ p, \Om,  w, r, \ti u,   $ be as in
Lemma 2.1.  Then 
there exists $ \al 
= \al ( p )  \in (0, 1) $ such that  $ \ti u $ has a  H\"{o}lder 
$ \al $ continuous representative in  $  B ( w,    r )  $ (also denoted $ \ti u$).   
 Moreover  if $ x, y \in B ( w, r )
$  then   
\[  | \ti u ( x ) - \ti u ( y ) | \leq \, c \,  ( | x - y |/ r )^\al  \, 
\max_{ B ( w, 2 r) } \,  \ti u. \]  }

 \noindent {\bf Lemma 2.3. } {\em Let $p, \Om,  w, r, \ti u,$ be as in Lemma 2.1  and let  $ \ti \mu $  be the measure
associated with $ \ti u $ as in (1.2). 
Then there exists $ c  $
 such that   
\[ c^{ - 1 } \, r^{ p - 2 } \, \ti \mu [ B ( w,  r/2 ) ]  \, 
\leq \, { \ds   \max_{B ( w,  r ) } \, \ti u^{ p - 1}   } \,  
 \, \leq \, 
 c \, r^{ p - 2 } \, \ti  \mu [ B ( w,  2 r  ) ].  \, \] } 
\

Using Lemma  2.3 we prove, \\

\noindent {\bf Lemma 2.4.} {\em Fix $ p, 1 < p < \infty, $ and let $
\hat u $ be  the   positive $p $ harmonic function
in Theorem 1.  Also, 
let $ u $ be the $p$ capacitary function for $ D = \Om \sem \bar B (
z_0, d ( z_0, \ar \Om )/2 ), $   defined  below
Corollary 1,  and let  $ \mu,  \hat \mu, $ be the measures corresponding 
to $ u, \hat u,  $ respectively.   Then  $  \mu , \hat  \mu $ are mutually
absolutely continuous. In particular,  Theorem 1  is  valid for $ \hat \mu $ if
and  only if  it is  valid for $  \mu. $}    \\  

\noindent {\bf Proof:}   We note that if $ \nu \not \equiv 0  $ 
is a   finite    
Borel measure on $ \ce $   with compact support, 
  then   
\eqn{2.5} \beq  \nu ( \ce \sem \Ga ) = 0  \mbox{ where }  \Ga = 
\left \{ z \in \mbox{ supp }\nu  :  \liminf_{ t \rar 0 }   \frac{ \nu ( B ( z,100 t
) ) }{ \nu ( B ( z, t ) )}  \leq  10^9  \right\}  
  \eeq   Indeed otherwise,  there exists a Borel set $ \La
\subset \ce  $ with $  \nu ( \La ) > 0 $ and the property that  
 if $ z \in  \La,  $  then  there exists $
t_0 ( z ) > 0 $  for which  
 \eqn{2.6}  \beq  \nu ( B ( z,  t ) ) \leq 10^{ -  8} \nu ( B ( z,
 100 t ) ) \mbox{ for }  0 < t <  t_0 ( z ). \eeq   
Iterating  (2.6)  it follows
that    \eqn{2.7} \beq  \lim_{t \rar 0 } \frac{ \nu ( B ( z, t )
)  }{
t^3 } = 0  \mbox{  whenever  $ z \in  \La.$}  \eeq
   Since    $  H^3 ( \ce ) = 0, $  we  deduce  from (2.7) 
  that    $ \nu ( \La ) =   0, $ which is a contradiction. 
 Thus  (2.5) is true.  
 
Now suppose that $ \mu,  \hat \mu $ are as in Lemma 2.4.  Let $
N_1 $ be a neighborhood of $ \ar \Om $ with  
\[   \ar \Om \subset  N_1 \subset \bar N_1 \subset  N. \] 
 Then from compactness and continuity of $ \hat u,
u,  $ there exists $  \hat  M   <  \infty  $ such that 
\eqn{2.8}  \beq  u \leq \hat M  \hat u \leq \hat M^2 u  \eeq 
 on $  \Om \cap \ar N_1.  $   
 From  (2.8)   and the boundary maximum principle for $ p $
harmonic functions we conclude that (2.8) holds  in $ \Om \cap
 N_1. $   In view of  (2.8)  and  Lemma  2.3  we see there exists $
 \hat r > 0, $   and a constant $ b  <  \infty, $ such that 
\eqn{2.9} \beq  \mu ( B ( w, s ) ) \leq b \hat \mu ( B ( w, 2 s )
)  \leq  b^2    \mu ( B ( w, 4 s ) )  \eeq 
whenever $ w \in \ar  \Om  $ and  $  0 < s \leq \hat r. $     
 We also note from  Lemma 2.3  that  supp $ \mu $ = supp $ \hat
\mu = \ar \Om. $ 
  
The proof of  Lemma 2.4 is by contradiction. 
  Let $   E \subset \ar \Om  $ be  a Borel set with $
\hat \mu ( E ) = 0.$    
If $ \mu ( E ) > 0, $ then from   properties of Borel
measures, and with $\Gamma$ as in  (2.5) with $ \nu = \mu, $     we
see  there exists a compact set $ K $ 
with  
\eqn{2.10} \beq   K \subset  E \cap \Ga  \mbox{ and } \mu ( K ) > 0.  \eeq 
 Given $ \ep > 0 $ there also exists an open set    
 $ O $   
with     
\eqn{2.11} \beq  E \subset O \mbox{ and }   \hat \mu ( O ) < \ep.   \eeq 
Moreover,  
we may suppose  
for each $ z \in K $  that there is  a  $  \rho = \rho ( z ) $
with   $ 
  0 < \rho ( z )  <
 \hat r / 1000,   $  $  \bar B ( z,  100 \rho (z)  ) \subset O , $ and  
\eqn{2.12}  \beq   \mu ( B ( z,  100 \rho  ) ) \leq 10^{10}   \mu
( B ( z,
   \rho  ) ).  \eeq  
 Applying  Vitali's covering theorem we then get   $ \{ B (
z_i,   r_i ) \}  $  with   $ z_i \in \ar \Om, 0 < 100 r_i  < 
\hat r $ and the property that  
\eqn{2.13} \begin{eqnarray}
 (a)&&\ \mbox{(2.12)  holds with $  \rho = r_i
$ for  each $ i,  $  }\notag \\  
(b)&&\ K \subset {\ds  \bigcup_{i } B ( z_i, 100 r_i ) }
\subset O,\notag\\
(c)&&\ B ( z_i, 10 r_i )  \cap  B ( z_j, 10 r_j )
 = \es \mbox{ when $ i \not = j $. } 
\end{eqnarray} 

Using   (2.9) and  (2.11) - (2.13),    it follows that 
\eqn{2.14} \begin{eqnarray}  
 \mu ( K ) &\leq&   \mu [  
\cup_{i } B ( z_i, 100 r_i ) ]
 \leq  \sum_{i }  \mu [ B ( z_i, 100 r_i ) ]  \leq 10^{10} 
\sum_{i }  \mu [ B ( z_i,  r_i ) ]\notag\\
 & \leq& 10^{10} 
 b \sum_{i }   \hat \mu [ B ( z_i, 10 r_i ) ]   \leq   10^{10}  \, b
\, \hat \mu ( O
) \leq 10^{10} \, b \,  \ep.  
\end{eqnarray}  Since $ \ep $ is arbitrary  we
conclude  that 
$ \mu ( K ) = 0,  $ which  contradicts  (2.10). Thus  $ \mu $ is
absolutely continuous with respect to $ \hat \mu. $ Interchanging
the roles of $ \mu, \hat \mu $ we also get  that $ \hat \mu $ is
absolutely continuous with respect to $ \mu. $   Thus Lemma 2.4
is true.  $ \Box $  
\\

\noindent   \section{ Proof of Theorem 1 (assuming Theorem 2).} From
Lemma 2.4 we  see that it suffices to prove Theorem 1 with 
$ \hat u, \hat \mu, $ replaced by $ u, \mu. $     
In this section we  prove Theorem 1  for $u$ under the
assumption that Theorem 2 is correct.   
Given Theorem 2  we can follow closely the argument in [BL05] from 
(6.9) on. However, our argument is    
 necessarily  somewhat  more complicated, 
 as in [BL05] we used the fact that  
 $ \mu $   was a doubling measure, which is not  necessarily true 
 when $ \Om $ is
simply connected.  We   claim that it suffices to prove Theorem 1
when    
\eqn{3.1}  \beq  z_0 = 0   \mbox{ and }  d ( z_0, \ar \Om ) = 2.
\eeq     Indeed, 
  let  $ \tau = d ( z_0, \ar \Om )/2 $ and put  $ T ( z ) =  z_0
+   \tau  z. $   If    
$ u' ( z ) = u ( T ( z ) ) $ for $ z \in D,  $ 
 then since the  $ p $ Laplacian is invariant under
translations,  rotations, dilations,  it follows that $ u' $ is
$ p $ harmonic in $  T^{-1} ( D ). $ Let $ \mu' $ be the measure
corresponding to $  u'. $   Then from (1.2)    it follows easily that 
\[  \mu' (  E )  =     \tau^{ p - 2}   \mu (  T ( E ) )  \mbox{ whenever $ E
\subset \rn{n}  $ is  a Borel set. } \] 
 This equality  clearly implies that  $ \hd \mu' = \hd \mu. $  
Thus we may assume that (3.1) holds. 
  Then $ B ( 0, 2 ) \subset \Om $
and  $ D = \Om \sem \bar B ( 0, 1). $   

Using Theorem 2 we have, for some $ c = c ( p ) \geq 1 $, that 
 \eqn{3.2} \beq 
c^{-1}\frac{u (z)}{d (z,\ar \Om)}\leq |\nabla u (z)|\leq
c \frac{u (z)}{d (z,\ar \Om)}  \mbox{ whenever $ z \in D.$ } 
 \eeq  
 Next set      
\[  v ( x ) =   \left\{ \bea{l} \max (  \log | \nabla u ( x ) |, 0 )
 \mbox {  when $ 1 < p < 2 $ } \\ 
  \max ( -  \log | \nabla u ( x ) |, 0 )  \mbox{  when  
 $ p >2. $ }  \ea  \right.  \]  
Then in [BL05] it is shown that    
\eqn{3.3} \beq  \int_{ \{ x: u ( x )  = t \} } \, | \nabla u |^{ p - 1}   \, 
 \exp \left[  \frac{w^2 }{ 2 c_+  \log (1 /t)} \right]  \,  d  H^1
x \, \leq \,  2 \,  c_+   \eeq 
 for some $ c^+ \geq 1. $   
 In [BL05],  $ c^+ $  depends on $ k, p,  $ but  only   because 
   the constant in (3.2) depends on $ k, p. $   So,  given Theorem
2,    $ c^+   = c^+ ( p ) $ in (3.3).   Next let  
   \[  \bea{c}   \xi ( t )  =  2 \sqrt{  \,  c_+ \, \log ( 1/t)  
 \,  \log \log (1/t ) }   \mbox{ for $ 0 < t < 10^{ - 6}, $
 } \\ \\     
 F ( t ) =  \{ x : u ( x ) = t \mbox{ and } v ( x )  \,    \geq
\, \xi  ( t )  \}. \ea  \] Then from (3.3) and 
weak type estimates we deduce
\eqn{3.4} \beq  \,  
 \int_{F ( t )} \, | \nabla u |^{ p - 1}  \, d H^1 x \, 
  \leq   2 c_+   \, [  \log ( 1/ t ) ]^{ - 2 }. 
    \eeq 

Next for $ A  $ fixed  with $ | A | $   large, 
we  define $ \la $   as in Theorem 1. Let $ a =  \frac{|A|}{
 2  \sqrt{c^+} } $  and  note that   
 \eqn{3.5} \beq   \la ( r )  =   \left\{  \bea{l}  r \, e^{   a  \xi  ( r ) } 
\, \mbox{ when } 1 < p <  2,  \\ 
r \, e^{ - a \xi  ( r ) } \mbox{ when } p > 2. \ea \right. \eeq 
To prove Theorem 1 when either $ 1 < p < 2 $ or $ p > 2, $ 
we intially allow $ a $ to vary but will later  fix   it as a
constant depending only on $ p,$ 
  satisfying several conditions.  
 Fix $ p,  1 < p <  2, $ 
 and let $  K \subset \ar \Om $ be a Borel set  with
$ H^\la   ( K ) = 0. $ Let $ K_1 $ be the subset of all  $ z \in  K $ with   
\[    \limsup_{r \rar 0}  \frac{ \mu ( B ( z, r ) )}{
\la ( r ) } 
 < \infty .  \] 
Then from the definition of $ \la  $ and   a covering argument
(see [Mat95, sec 6.9]),  it is easily shown that $   \mu  ( K_1 ) 
= 0. $  
Thus to prove   $ \mu ( K ) = 0, $ it suffices to show 
 $ \mu ( E ) = 0 $ when $ E $ is Borel and is equal  $ \mu $ almost
everywhere to  the set of all points in 
$ \ar \Om $ for which   
\eqn{3.6} \beq   \limsup_{r \rar 0}  \frac{ \mu ( B ( z, r ) )}{
\la ( r ) } 
 = \infty .  \eeq 
 Let  $ G  $ be the set of all $ z $  where $ (3.6) $ holds.  
Given 
$  0 < r_0  < 10^{ - 100},  $   we  first show for each $ z \in G $ that
there exists   $ s = s ( z ), 0  < s/100 < r_0, $ such that 
\eqn{3.7}  \beq  \mu (  B ( z,  100 s) ) \leq 10^9  \mu ( B ( z, s ) ) 
\mbox{ and }   \la (100  s ) \leq     \mu ( B ( z, s ) ) .  \eeq    
In fact let $ s \in ( 0, r_0 ) $ be the first point  starting
from $ r_0 $  where  
\[  \frac{ \mu (  B ( z,  s) ) }{ \la ( s ) } \geq 10^{20}    
  \min  \left\{ \frac{ \mu ( B ( z, r_0  ) )}{ \la ( r_0 ) }, \, 
1   \right\} . \]
 From (3.6) we see that $ s $  exists.   Using $ \la ( 100  r ) \leq
 200  \la ( r ), 0 < r < r_0/100, $   it is also easily checked that  
  (3.7) holds.   From (3.7) and  Vitali again,      
we get  
  $ \{ B (
z_i,   r_i ) \}  $  with   $ z_i \in  G,  0 < 100 r_i  < 
  r_0,   $ and the property that  
\eqn{3.8} \begin{eqnarray} 
(a)&&\ \mbox{(3.7)  holds with $
 z = z_i, s  = r_i, 
$ for  each $ i,  $}\notag \\ 
(b)&&\ G \subset {\ds  {\ds \bigcup_{i } }  B ( z_i, 100 r_i ) }\notag\\
(c)&&\ B ( z_i, 10 r_i )  \cap  B ( z_j, 10 r_j )
 = \es \mbox{ when $ i \not = j $. }
 \end{eqnarray}  
 Let  $ t_m = 2^{ - m } $ for $ m =   1, 2,  \dots.  $   Given $
i, $ we claim  there exists   $ w_i \in   B ( z_i,
5  r_i )  $  and  $ m = m ( i ) $ with       
\eqn{3.9 } \begin{eqnarray}
 (\al)&&\ u ( w_i ) = t_m  
\mbox{ and } d ( w_i, \ar \Om ) \approx  r_i\notag \\ 
 (\be)&&\ \mu [ B ( z_i,  10 r_i ) ]/ r_i  \approx    [ u ( w_i ) / d ( w_i, \ar
\Om ) ]^{ p - 1} \approx   | \nabla u ( w ) |^{ p - 1}\notag\\  
&& \mbox{ whenever $ w \in B ( w_i, d ( w_i, \ar
\Om )/2 ) .$ } 
\end{eqnarray}
   In (3.9) all proportionality
constants depend only on $p.$  To prove (3.9)   
choose $ \ze_i  \in \ar B ( z_i, 2  r_i ) $ with    
  $ u ( \ze_i )  =  {\ds \max_{ \bar B (  z_i,  2   r_i ) } }   u. $   Then
 $ d ( \ze_i, \ar \Om ) \approx  r_i, $   
 since  otherwise,  it would
follow  from Lemma 2.2 that  $ u ( \ze_i ) $ is  small in comparison to 
  $ {\ds \max_{ \bar B (  z_i,  5   r_i ) } }  u . $   However  
 from (3.8) $ (a) $ and Lemma 2.3,  these two
  maximums  are proportional with  constants
depending only on $ p. $    Thus  $ d ( \ze_i,  \ar \Om  )  
\approx r_i. $         
Using this fact,  (3.2), (3.8) $(a),$ and Lemma
2.3, once again we get (3.9) $ (\be) $ with $ w_i $ replaced by 
$ \ze_i. $  If   $ t_m  \leq  u ( \ze_i  )  <  t_{m-1} $ we let $ w_i $ be
the first point on the line segment  connecting $ \ze_i $ to a
point in $ \ar \Om \cap \ar B (  \ze_i, d ( \ze_i, \ar \Om ) ) $
where $ u  = t_m. $   From our construction, 
 Harnack's inequality,  and Lemma
2.2  we see that   (3.9) is true.   

Using (3.8), (3.9),    we  deduce  for $ 1 < p < 2 $ that 
\eqn{3.10}  \beq  v ( z ) =   \log  | \nabla u ( z ) | \geq  \, 
a \, \xi  ( 100  r_i ) / \ti c  
\mbox{ on }  B ( w_i,  d ( w_i, \ar \Om )/2 ) \eeq  
 where $ a$ is as in (3.5).   Next we note  that 
\eqn{3.11} \beq      H^1 [ B ( w_i, d ( w_i, \ar \Om )/2 )  \cap
\{ z : u ( z )
= t_m \} ]   \geq d ( w_i, \ar \Om )/2  \eeq  as we see from the maximum
principle for $p $ harmonic functions,  a connectivity
argument and basic geometry.    Also,     we can
use  (3.8) $(a)$  to estimate $ t_m  $ below in terms of  $ r_i  $
 and     Lemma 2.2   to  estimate $ t_m $ above in terms of $
r_i. $  Doing this 
we find
for some $ \be = \be ( p ), 0 < \be < 1,  \bar c =  \bar c ( p ), $ that 
 \eqn{3.12} \beq   r_i \leq   \bar c  \,  t_m^\be \leq \bar c^2
\,  r_i^{\be^2} \, .
\eeq 
Using  (3.8)-(3.12)     
 we 
conclude, for $ a $ large enough, that  
\eqn{3.13} \beq  
\mu [ B ( z_i, 10 r_i ) ] 
\, \leq \,  c \, \int_{  F ( t_m ) \cap B ( z_i, 10 r_i ) }  
 | \nabla u
|^{ p - 1} \, d H^1 . \eeq 
Using  (3.8),    
 (3.12), (3.13), and (3.4)    it follows for $ c $ large enough that 
\eqn{3.14} \begin{eqnarray}
  \mu ( G  )  &\leq& {\ds  \mu \left( \bigcup_i B (
z_i,  100 r_i )  \right) }   
 \leq  \, { \ds  10^9 \,   \sum _{i}    
  \, \mu  [  B ( z_i, 10 r_i ) ] }\notag\\
    & \leq & c { \ds \sum_{ m = m_0 }^\infty 
 \int_{  F ( t_m )  }  
   | \nabla u |^{ p -  1} d H^1 x \, }  
  \leq  \, c^2   { \ds \sum_{m = m_0}^\infty } m^{- 2} \, \leq  
c^3  m_0^{ - 1}
\end{eqnarray}
where $   2^{ - m_0 \be } =   \bar c  \, r_0^{\be^2}. $ Since $ r_0 $ can be
arbitrarily small we see from (3.14)  that $ \mu ( G ) = 0. $
This equality and the  remark above  (3.6) yield $ \mu ( K ) = 0.
$ Hence $ \mu $ is absolutely continuous with respect to $ H^\la  $
and Theorem 1 is true for $ 1 < p < 2. $ 

Finally to prove Theorem 1  for  $ p > 2, $ we show 
there exists  a Borel set  $  \hat K \subset  \ar \Om $ 
 such that  
\eqn{3.15} \beq  \mu ( \hat K ) = \mu ( \ar \Om )  \mbox{ and }  
\hat K   \mbox{ has $ \si $ finite $ H^\la $ measure.}  
\eeq  
  In fact let $ \hat K $ be the set  of all  $ z \in   \ar \Om  $ with   
\eqn{3.16} \beq   \limsup_{r \rar 0}  \frac{ \mu ( B ( z, r ) )}{
\la ( r ) } 
 >  0.  \eeq  
  Let   $ \hat K_n $ be the subset of $ \hat K $ where the above 
$ \limsup $ is  greater than $ 1/n. $ Then  from  the definition of $ \la $  
 and a Vitali covering
type argument (see [Mat95, ch 2])  it  follows  easily  that
  \[ H^\la ( \hat K_n ) \leq  100 n  \mu ( \hat K_n). \]   
  Since  $ \hat K =  \cup_n \, \hat K_n   $ we conclude that 
$ \hat K $  is $ \si $ finite with respect to $ H^\la $ measure.
Thus   to prove (3.15) it suffices to show 
 $ \mu  (  \hat G  ) = 0 $ where  $ \hat G  $ is   equal  
 to  the set of all points in 
$ \ar \Om $ for which  
\eqn{3.17}   
\beq   \lim_{r \rar 0}  \frac{ \mu ( B ( z, r ) )}{
\la ( r ) } 
 = 0.  \eeq 
Given 
$  0 < r_0  < 10^{ - 100}  $   we   argue as in  the proof of
(2.5) to  deduce for each $ z \in \hat G $ the 
 existence of    $ s = s ( z ), 0  < s/100 < r_0, $ such that 
\eqn{3.18}  \beq  \mu (  B ( z,  100 s) ) \leq 10^9  \mu ( B ( z, s ) ) 
\mbox{ and }  \la (  s ) \geq   \mu ( B ( z, 100 s ) ) .  \eeq    
    Using (3.18) and  once again applying Vitali's covering lemma   
we get  
  $ \{ B (
z_i,   r_i ) \}  $  with   $ z_i \in  \hat G,  0 < 100 r_i  < 
  r_0,  $ and the property that  
\eqn{3.19}\begin{eqnarray} 
(a)&&\ \mbox{(3.18)  holds with $
 z = z_i, s  = r_i
$ for  each $ i,  $  } \notag\\ 
(b)&&\ \hat G \subset {\ds  {\ds \bigcup_{i } }  B (
z_i, 100 r_i ),  }\notag\\ 
(c)&&\ B ( z_i, 10 r_i )  \cap  B ( z_j, 10 r_j )
 = \es \mbox{ when $ i \not = j $. } 
 \end{eqnarray}   
Let  $ \He   $ be the set of all indexes, $ i, $ for which 
$ \mu ( B ( z_i, 100 r_i ) ) \geq r_i^3 $ and let $ \He_1 $ be
the  indexes for which this inequality is false.  Arguing as in
(3.14) we obtain  
 \eqn{3.20}  \begin{eqnarray}  
\mu (  \hat G )  \leq    \bigcup_{i \in \He }   \mu ( B ( z_i, 100 r_i ) )
 +  \sum_{i \in \He_1 } r_i^3 \, \leq  
  10^9 \, \bigcup_{i \in \He }   \mu ( B ( z_i, 10 r_i ) ) \, + \, 100
r_0 \, ( H^2 ( \Om ) + 1 ). \end{eqnarray} 
If $ i \in \He, $  we   can repeat the argument  after (3.8) to
get  (3.9).  (3.9) and (3.8) $(a)$ imply  (3.10) for 
 $ w  = - \log | \nabla u |. $      Also since $ i  \in \He $ we can use
(3.9) to estimate $ t_m $ from below  in terms of $ r_i $  and
once again use Lemma 2.2  to  estimate $ t_m $  from above in terms of
$ r_i. $    Thus  (3.12) also holds for some $ \be,  \bar c  $
depending only on $p. $  (3.10) - (3.12)  imply (3.13) 
for $ a $ (as in (3.5)) suitably large.  
 In view of   (3.20), (3.13),  and (3.4) we have  
\eqn{3.21} \begin{eqnarray}  \mu ( \hat G )  - 100 r_0 ( H^2 ( \Om ) + 1 )  
 &  \leq& {\ds  \mu \left( \bigcup_i B ( z_i, 100 r_i )  \right) }   
   \leq  \, { \ds  10^9 \,   \sum _{i \in \He}    
  \, \mu  [  B ( z_i, 10 r_i ) ] } \notag\\
   & \leq& c { \ds \sum_{ m = m_0 }^\infty 
 \int_{  F ( t_m )  }  
   | \nabla u |^{ p -  1} d H^1 x \, }  
  \leq  \, c^2   { \ds \sum_{m = m_0}^\infty } m^{- 2} \, \leq  
c^3  m_0^{ - 1}  
\end{eqnarray}  
where $   2^{ - m_0 \be } =   \bar c  \, r_0^{\be^2}. $ Since $ r_0 $ can be
arbitrarily small we conclude  first  from (3.21)  that $ \mu ( \hat G ) = 0 $
 and thereupon that (3.15)  is valid.    Hence $ \mu $ is concentrated on a  set of $  \si $ finite 
 $ H^\la  $  measure when  $  p > 2. $ The proof of Theorem 1
is now complete  given  that Theorem 2 is true.  $ \Box $

\section{ Preliminary Reductions for  Theorem 2.}
  Let $ u $ be the  $p$ capacitary function for 
$ D = \Om \sem B ( z_0, d ( z_0,  \ar \Om )/2 ). $  We extend $ u
$ to   $ \ce $   by putting $ u \equiv 1 $ on   
$  \bar B ( z_0,  d ( z_0, \ar \Om )/ 2 ) $ and $ u
\equiv 0 $ in  $ \ce \sem \Om. $       
 We shall need
some more basic properties of $u$.  Again references for proofs can be found in
 [BL05]. \\ 

\noindent  {\bf  Lemma 4.1.}
  { \em  If  $ z = x + i y, i = \sqrt{-1},  x, y  \in \re, $  then  
  $u_{z}= (1/2)
(u_{x}-iu_{y})$ is a 
quasi-regular mapping of $ D $ and $\log |\nabla u|$
is a weak solution to a linear elliptic PDE in divergence
form in $D$.  Moreover, positive weak solutions to this PDE in  $ B ( \ze, r )
 \subset D $ satisfy the Harnack
inequality 
\[
\max_{B( \ze,r/2)}h\leq \tilde{c}\min_{B ( \ze ,r/2)}h
\] 
where $\tilde{c}$ depends only on $p$. } 
\\

 \noindent {\bf Lemma 4.2.} 
  {\em $u$ is real-analytic in $ D, $ $\nabla u\neq 0$ in
$ D,  $ and $\nabla u$ has a H\"older continuous
extension to a neighborhood of $ \ar B ( z_0, d ( z_0, \ar \Om)/2).$
Moreover, there are constants $ \be, 0 < \be  <  1,  $ and $  \hat{c}\geq 1$, depending only on $p$,
such that
\[
|\nabla u (z) - \nabla u ( w) | \leq  
\hat{c} 
 \left( \frac{ | z - w | }{d ( z, \ar \Om )} \right)^\be \,   
\max_{B ( z,  d( z, \ar \Om )/2 )}  \,  | \nabla u |  \,   
 \leq   
  \,  \hat{c}^2  \, 
 \left( \frac{ | z - w | }{d ( z, \ar \Om )} \right)^\be \,   
  \frac{u (z)}{d (z, \ar \Om)} \]
whenever $ w \in   D \cap  B (  z, d ( z, \ar \Om )/2 ).  $  
Finally    
\[
 \hat c \, |\nabla u ( w) | \geq   
 \,  \frac{u (w)}{d ( w, \ar \Om)}
 \mbox{ for  $ w \in   D \cap  B (  z_0,   3 d ( z_0, \ar \Om
)/ 4 ).  $}  
\] } \

Using Lemma 4.2 we see that Theorem  2 is true when 
$ z \in  D \cap B ( z_0,  3 d ( z_0, \ar \Om )/4 ). $ Thus it is enough to 
prove Theorem 2 with $z=z_1$ for 
\eqn{4.3} \beq  z_1  \in 
   D \sem B ( z_0,  3 d ( z_0, \ar \Om )/4 ). \eeq 
Recall  the  definition of the hyperbolic distance $\rho_\Om$  
for a simply connected domain $ \Om$ (see [GM05]). 
Then $ \rho_\Om(z_{1},z_{2})$, $z_{1},z_{2}\in\Omega$, is comparable to the
quasi-hyperbolic distance 
\[ Q_{\Om} (z_{1},z_{2}) : = \inf\int_\ga \frac{|dz|}{ d (z,\ar \Om)}
 \]
where the infimum is taken over all the paths $\ga \subset\Om$
connecting $z_{1}$ to $z_{2}$. More specifically,  
\eqn{4.4}  \beq \rho_{\Om}\leq Q_{\Om}\leq
4\rho_{\Om} \eeq  as follows from the Koebe estimates
 \eqn{4.5} \beq 
\frac{1}{4}|f^{\prime} (z)| (1-|z|^{2})\leq  d (f (z),\ar \Om)\leq
|f^{\prime} (z)| (1-|z|^{2}),  \, z \in B ( 0, 1 ),   
\eeq 
whenever $f: B ( 0, 1) \rightarrow \Om$ is a conformal map, (see
Theorem I.4.3 in
  [GM05]). In the following we will often use the following distortion
estimate, which also follows
from Koebe's Theorem, (see (I.4.17) in [GM05]), for conformal maps $f: B ( 0, 1)  \rightarrow \ce $. 
For $z_1,z_2\in D$,
 \eqn{4.6} \begin{equation}
\rho_{\Om }(z_1,z_2)\leq A_{1}\Longrightarrow|f^{\prime}
(f^{-1}(z_2))|\leq A_{2}|f^{\prime} (f^{-1}(z_1))|
\end{equation}
for some constant $A_{2}$
depending only on $A_{1}$. Note also that (4.6) implies that $d(z_2,\partial\Omega)\leq A_3d(z_1,\partial\Omega)$
for some constant $A_{3}$
depending only on $A_{2}$. The same holds if  $ f $ is a conformal
mapping of the  
upper half-plane $\Hh $.
Our  main lemma  in the proof of Theorem 2 is the following. \\

\noindent  {\bf Lemma 4.7.}
{\em There is a constant $ C$, depending only on $p$, 
such that if $z_{1}$ is as in (4.3) then there exists $z^\star\in \Omega$
with $u (z^\star)=u(z_{1})/2$ and $\rho_\Om (z_{1},z^\star)\leq
C$.}  \\

    Assuming for  the  moment that  Lemma 4.7    is proved we get
Theorem 2 from the following argument.  
 Let  $  \Ga $ be the
hyperbolic geodesic connecting $ z_1 $ to $ z^*. $ If 
$ \Ga \cap B ( z_0, 5 d ( z_0, \ar \Om )/8 ) = \es, $ 
we put 
$ \ga = \Ga. $  Otherwise,   $   \ga = \ga_1 + \ga_2 + \ga_3  $
where $ \ga_1 $ is the subarc of $ \Ga $ joining $ z_1 $  to the
first point, $  P_1, $  where $ \Ga $ intersects   
$ \ar B ( z_0, 5 d ( z_0, \ar \Om )/8 ) $;  $ \ga_2 $ is  the
short  arc
 of  $ \ar B ( z_0, 5 d ( z_0, \ar \Om )/8 ) $  joining $ P $ to
 the last point, $ P_2, $  where $ \ga $ intersects 
$ \ar B ( z_0, 5 d ( z_0, \ar \Om )/8 )$; and   finally $ \ga_3 $
joins  $ P_2  $ to $ z^* . $
  Using  (4.3)-(4.6), one
sees that  
  \eqn{4.8} \beq    H^1 ( \ga ) \leq c d (z_{1},\ar \Om)  
\mbox{ and }  
 d (\ga ,\ar \Om)\geq
c^{-1} d (z_{1},\ar \Om), \eeq  where $ c = c ( p ). $  
Thus 
\[
\frac{1}{2}u (z_{1})\leq u (z_{1})-u (z^{\star})\leq \int_{\ga}|\nabla
u (z)||dz|\leq c H^{1} (\ga) \, \max_{\ga}|\nabla u|\leq C d
(z_{1},\ar \Om)\max_{\ga}|\nabla u|.
\]
So for some $\zeta\in \ga$,
\eqn{4.9} \beq 
c^{\star}|\nabla u (\zeta) |  \, \geq  \, \frac{u (z_{1})}{d (z_{1},\ar \Om)}
\eeq 
 where  $c^{\star}\geq 1$ depends only on $p$.
 Also  from (4.8) 
 we deduce the existence of balls $\{B ( w_j, r_j \}_{j=1}^{N} , $  with 
 $ w_j \in \ga  $ and 
 \eqn{4.10} \begin{eqnarray} (a)&&\ B ( w_j, r_j/4 ) \cap B ( w_{j+1} , r_{j+1}/4  )  \not = \es 
\mbox{ for $ 1 \leq j \leq N - 1, $ }\notag  \\ 
(b)&&\  r_j  \,  \approx  \,  d( B ( w_j, r_j ) ,\ar \Om) \,  
\approx  d ( z_1, \ar \Om ),\notag  \\
(c)&&\ \ga \subset  \bigcup_j B ( w_j, r_j/4 ),
\end{eqnarray}  
where  $ N$  and proportionality constants depend only on $p$.
 Observe from (4.10) and  Harnack's inequality applied to $ u $
(see Lemma 2.1)  that  $ u ( z ) \approx u ( z_1 ) $  when 
$ z \in \cup_j B ( w_j, r_j ). $  In view of Lemma 4.2, (4.10),  it follows
 for some $ c = c ( p ) $ that  
\eqn{4.11} \beq   | \nabla u ( z )  | \leq c u ( z_1 )/ d ( z_1, \ar
\Om )  \mbox{ when } z \in \bigcup_{j} B ( w_j, r_j/2 ) \, . \eeq 
 From (4.11)  we see that if $ c = c ( p ) \geq 1 $ is large
enough and     
\[
h (z) = : \log \left(\frac{c \,  u (z_{1})}{d
(z_{1},\ar\Om) \,  |\nabla u (z)|} \right) \mbox{ for $ {\ds  z \in  
\bigcup_j B ( 
w_i, r_i/2 ) } $ }  
\]
 then    $ h  >  0 $ in 
$  \cup_i B ( w_i, r_i/2 ). $ 
Choose $ i, 1 \leq i \leq N, $ so that    $ \zeta \in  B ( w_i,
r_i/4  ). $  
Using (4.9) we have  $  h ( \ze ) \leq c. $ Applying 
 the Harnack inequality in Lemma  4.1  to $ h $  in  
$ B ( w_i, r_i/2 ) $  
 we get 
\eqn{4.12}  \beq c   | \nabla u |  \geq  u ( z_1 )/d ( z_1, \ar \Om)   
  \mbox{   in    $ B ( w_i, r_i/4).  $  } \eeq   From 
 (4.10) we see that the argument leading to (4.12)  can be
repeated in a chain of  balls connecting  $ \ze $ to  $ z_1. $
Doing this and using $ N = N ( p ), $  we get Theorem 2. $ \Box $  \\ 

  In the proof of Lemma  4.7 we may  assume without loss of generality that $\ar \Om$ is an analytic Jordan
curve, as the constant in this lemma  will depend only on $p$. 
Indeed,  we can approximate $\Om$  
 by an increasing sequence
of analytic Jordan domains $\Om_n \subset \Om $, and   apply
Lemma 4.7 to
$ u_n $ the $p$ capacitary function for $ D_n = \Om_n \sem B (
z_0, d ( z_0, \ar \Om )/ 2  ). $ 
 Doing this  and letting $ n \rar \infty, $    we get 
Lemma 4.7 for $ u,  $ since  by  
Lemmas 2.2, 4.2,  there are  
subsequences  of $u_n,$
$\nabla u_n,$  converging  to $u,  \nabla u,$ respectively,
uniformly on compact subsets of $\Om$.

\subsection{Outline of the proof  of Lemma 4.7.}
  To prove  Lemma 4.7   
It will be useful to transfer the problem to the upper half-plane
$\Hh $ via the
Riemann map $f:\Hh \rightarrow \Om$ such that $f (i)=z_{0}$ and
$f(a)=z_1$ where $a=is$ for some $0<s<1$.  We note that $ f $  
 has a continuous extension to $ \bar \Hh, $ since $   \ar \Om  $
is a Jordan curve.  We also let $U = u\circ f$,
and note that $ U $  satisfies a  maximum principle and Harnack's inequality. Consider the box
\[
Q (a) =  \{z=x+iy: |x |\leq  s,\ 0<y< s \}.
\] 
We will show that $Q (a)$ can be shifted to a nearby box $\tilde{Q}
(a)$ whose boundary in $\Hh $ we call $\xi$. 
It consists of the horizontal segment from $x_1+is$
to $x_2+is$, and the vertical segments
connecting $x_l+is$ to $x_l$ for $l=1,2$.  $x_1, $ $x_2,$ are
chosen  to satisfy  $- s <x_1 <- s/2, \,  
 s/2<x_2< s. $ Let $ f ( x_j ) = w_j, j = 1, 2. $  $ \ti Q ( a ) $ will be constructed to have several nice properties. In particular,  we  will prove that 
$U\leq A U(a),$
on $\xi$, and hence, by the maximum principle, $U\leq A U(a)$ on $\tilde{Q}(a)$, for
some constant $A$ depending only on $p$. In other words,
if we let $\sigma\defeq f(\xi)$ and $\Om_1\defeq f(\tilde{Q}(a))$, then we will prove that
 \eqn{4.13} \begin{equation}
\label{eq:maxpr2}
u\leq Au(z_1)
\end{equation}
on $\si$ and hence in $\Om_1$. Moreover, we will prove that
 \eqn{4.14} \begin{equation} \label{eq:silength}
H^1 (\sigma) \leq C_1 d ( z_1,\ar \Omega)
\end{equation}
for some absolute constant $C_1$ depending only on $p$. Furthermore,we will establish the existence of $w_0=f(x_0)$, for some
$|x_0|<s/4$, such that $|w_0-z_1|\leq  C_2 d (z_1,\ar \Om)$ and such that
 \eqn{4.15} \begin{equation}\label{eq:w0}
d (w_0,\sigma)\geq d(z_1,\ar \Om)/C_2
\end{equation} where  $C_2  $ is an other absolute constant. In addition we will 
construct a Lipschitz curve $\tau:[0,1)\rightarrow
  \Om_1$ with $\tau(0)=z_1$ and $\tau(1)=w_0$, which
satisfies the cigar condition
 \eqn{4.16} \beq  
\min\{ H^1(\tau[0,t]), H^1(\tau[t,1])\}\leq
  C_3d(\tau(t),\ar \Om), 
\eeq 
for $0 \leq t\leq 1$ and some absolute constant $C_3.$ 

\noindent To briefly outline the construction of $\tau$ we note that we construct $ \tau $
 as the image under $f$ of a
polygonal path $$\lambda =  
{ \ds \sum_{k=1}^\infty   \la_k  } \subset \tilde{Q} (a),$$ starting at $a$ and
tending to $x_{0}$ non-tangentially. 
The segment $\lambda_k$,  $ k = 1, 2, \dots,  $ 
 joins $ a_{k-1} $ to $ a_k $ and consists of a horizontal line
segment followed by  a  downward
pointing vertical
segment. More precisely, fix $ \delta, $  
$0<\delta<10^{-1000}$ and put $  \de^* = e^{ - c^* / \de },  t_0
= 0, s_0 = s,   a_0 = t_0 + i s_0 = a. $  
In our construction we initially
allow $ \de $ to vary but shall  fix  $  \de $ in (5.3) to be a small
positive absolute constant  satisfying several conditions.  
 Also, $ c^*  \geq 1 $ is an absolute constant
which will be defined  in   Lemma  4.26.   
    Then $\lambda_1$ consists of 
 the horizontal segment from $a_0 $ to $t_1+is_0$ followed by the vertical
segment from $t_1+is_0$ to $a_1 =  t_1+i\de^*  s_0$.   Put
 $ s_1 =  \de^*  s_0. $   
Inductively,   if  $a_{k-1} =  t_{k-1}+is_{k-1}$ has been defined, 
 then $ \la_k $  consists of the horizontal line segment
joining  $ a_{k-1} $ to   
$t_k+is_{k - 1}, $ followed by the vertical line segment
connecting 
$t_k+is_{k - 1} $  to  $ a_k  = t_k + i s_k,  $ where 
$ s_k = \de^*  s_{k-1}.  $  
   Moreover the numbers
$t_k, k = 1, 2, \dots, $ are chosen  in such a way that
 \eqn{4.17} \begin{equation}\label{eq:key}
  | t_k - t_{k - 1} | \leq s_{k - 1}  
\mbox{ and }  \int_{0}^{s_k}   
 | f' ( t_k + i  \tau   ) |  d \tau  \,   \leq \, 
 \delta  \, d (f (a_{k-1}), \ar \Om).
\end{equation}
    Existence of  $ (t_k) $ will be shown in the paragraph after (5.2). 
  Letting $\tau_k=f(\lambda_k)$ and $z_k=f(a_{k-1} ), k = 1, 2,
\dots, $  we note that 
(4.17) and our construction imply     
 \eqn{4.18} \begin{equation}\label{eq:taudistbd}
d(z_{k+1},\ar \Om)\leq \delta d (z_{k},\ar \Om).
    \end{equation}
    For $w\in \lambda_k$, (4.6) and our construction give a
constant $\bar c$, depending only on $\delta$  and
$p$, such that \[ 
\bar c^{ - 1}  |f' ( a_{k-1}  ) |   \leq  |f^{\prime}(w)| \leq   
 \bar c   |f^\prime 
(a_{k - 1}) | \mbox{ whenever } w \in \la_k . \] Consequently for some constant  $  c\geq 1$, depending only on $\delta$  and
$p$,
 \eqn{4.19} \begin{equation}\label{eq:taulength}
  c  \, d ( w, \ar \Om )  \geq  d ( z_{k}, \ar \Om )
 \mbox{ when $ w \in \tau_k $ and }  H^1 (\tau_k)\leq c(\de)  d
(z_{k},\ar \Om)
\end{equation}
for $k=1,2,3,...$

Putting (4.18) and (4.19)  together we see
that if  $w  = \tau ( t ) \in \tau_k, $  then for some 
$ c_+ \geq 1,$ depending only on $\delta$  and
$p$, 
 \[   
|w -w_0|\leq   H^1 ( \tau [ t, 1] )    \leq   {\ds
\sum_{j=k}^{\infty} }  \,   H^1 (\tau_j) \,   \leq
 c_+     d ( w, \ar \Om )  \leq c^2_+ 
\delta^{k-1} d (z_{1},\ar \Om).
 \]
Using  this equality and (4.19)  we conclude that  $ \tau $ satisfies the cigar
condition in  (4.16) with a constant depending only on $ \de, p$.   
 
   To show the existence of $ z^* $  in Lemma 4.7,   
we  suppose    $ \de  > 0  $   is now  fixed as in (5.3)    
 and suppose that  
  $\lambda$ is parametrized by $[0,1]$ with
$\lambda(0)=a$ and $\lambda(1)=x_0.$ Let
\[
t^\star=\max\{t: U(\lambda(t))= {\ts \frac{1}{2}} U(a)\}
\]
 and put   $a^\star$ =  $\lambda(t^\star)$ and $z^\star=f(a^\star)$.
 If $ \rho  = d ( w_0, \si ), $ Then from   the definition of $ \Om_1 $
above  (4.13)  we have  
\[
B(w_0,\rho)\cap\Om\subset \Om_1.
\]
so  from  Lemma 2.2  
applied to the restriction of $u$ to $\Om_1,$ (4.13), 
  (4.15), and (4.16)   we  deduce  for some $ \ti c 
= \ti c ( p ) $ that    \[
 \frac 12u(z_1) = u(z^\star)\leq  \ti c  \left(\frac{d( z^*, \ar
\Om)}{\rho}\right)^{\alpha} \max_{B(w_0,\rho)\cap \Om}u \leq
 \ti c^2   A  \left(\frac{d(z^*,\ar \Om)}{d (z_{1},\ar \Om)}\right)^{\alpha}u(z_1). 
\]
 Thus  \[ d ( z_1, \ar \Om ) \leq  
   c  d ( z^* , \ar \Om )   \mbox{ for some  $ c = c ( p )$. }
  \]  
This inequality and  (4.16) imply that there is a  chain of
 $ N = N ( p ) $ balls  (as in (4.10))  connecting 
$ z_1 $ to $ z^*. $   Using  this implication and  once again (4.4)  we
conclude that  $ \rho_\Om  ( z^*, z_ 1 ) \leq c. $ This completes our
outline of the proof of  Lemma 4.7.

  To  finish the proof of  Lemma 4.7 we show there exists  $ \de
 > 0,  
\si, \tau, c^*,  (\tau_k)_1^\infty, $ 
 for which  
  (4.13) - (4.15) and (4.17),  are true. 

\subsection{  Several Lemmas.  } 
 To set the stage for the  proof of (4.13) - (4.15) and (4.17)  
we shall need several lemmas.  To this end  
define, for $b\in\Hh$, the interval $I (b)\defeq [\rea b-\ima
b,\rea b+\ima b]$.
\\

\noindent{\bf  Lemma 4.20.} 
{ \em  There is an absolute constant $ \hat C  $ such that if $f$ is
univalent on $\Hh$ and $b\in\Hh$, then } 
\[
\int\int_{\Hh}\frac{|f^\prime(w)|}{|f(w)-f(b)|}dA(w)\leq  \hat C
\ima b.
\]  
{\bf Proof of  Lemma  4.20:}  
The proof is left as an exercise. Hints are provided  in
  problem 21 on page 33 of [GM05], where the  case for
functions $g$ univalent on $ B ( 0, 1) $ with $ \rea g\neq 0$ is discussed. The same arguments
give the result for univalent functions  $ g $ on $ B ( 0, 1 ) $ 
with $g (0)=0$ and then  Lemma 
4.20  is obtained by applying the result to $g=f\circ
M_{b}$ where $M_{b} (z)=i \ima b (1+z)/ (1-z)+\rea b$.
$ \Box $ \\ 

\noindent {\bf  Lemma  4.21.}  
{\em There is a set $E (b)\subset I (b)$ such that for $x\in
E(b)$ } \eqn{4.22}  
  \begin{equation}\label{eq:length}
\int_0^{{\ts \ima b }}|f^\prime(x+iy)|dy\le C^\star d(f(b),\ar \Om)
\end{equation} 
for some absolute constant $C^\star,$ and also
\eqn{4.23} \begin{equation}\label{eq:h1i}  H^1 ( E ( b ) )  
 \geq  ( 1 -  10^{ - 100} ) H^1 ( I ( b ) ). \end{equation}

 Note that we could achieve Lemma 4.21 by invoking known results
in the literature, such as the result in [BB99]  related to previous
theorems of Beurling and Pommerenke (see [P75], Section 10.3).  
 For completeness we  give  an
alternative proof of Lemma 4.21  based on  Lemma 4.20.  
\\ \\  

 \noindent {\bf Proof of Lemma  4.21:}
Let $\ell$ be a large positive integer that will soon be fixed as an
absolute number and  let \[  T  = T(b) = 
\{ z = x + i y :  |  x  |
< \ima b:\  y = \ima b \}  \] be the top of the box 
 $ Q ( b ) $  defined  at the beginning of subsection 4.1.  Set 
\[
K=K(b)\defeq\{x\in I(b): |f(x+it)-f(b)|>2^\ell |f^\prime(b)|\ima
b\mbox{
  for some }0<t<\ima b\}.
\]
Note that
\[
|\partial_y \log|f(z)-f(b)||\leq \frac{|f^\prime(z)|}{|f(z)-f(b)|}.
\]
Also, for $z$ in the top $T$,
\[
|f(z)-f(b)|\leq 1000|f^{\prime}(b)|\ima b.
\]  
Thus,
\[
\int_0^{\ima {\ts b}}\frac{|f^\prime(x+iy)|}{|f(x+iy)-f(b)|}dy\geq
\frac{\ell }{C},
\]
whenever $x\in K$. Integrating both sides over $ K  $ and using
Lemma 4.20  we therefore find that
\eqn{4.24} \begin{equation}\label{eq:h1ka}
H^1(K)\leq C\frac{\ima b}{\ell}.
\end{equation}

%Consider $M(\chi_K)$, the Hardy-Littlewood maximal function of $K$,
%and put
%\[
%\hat{K}\defeq\left\{x\in \R: M(\chi_K)(x)>\frac{1}{\sqrt{\ell}}\right\}.
%\]
%Then, from weak-type estimates and (\ref{eq:h1ka}),
%\begin{equation}\label{eq:h1khat}
%H^1(\hat{K})\leq \tilde{C}\frac{\ima a}{\sqrt{\ell}}
%\end{equation}
%where $\tilde{C}$ is an absolute constant.
Next for we define a function $g(x)$ for $x\in I(b)$ as follows. If 
$x\in I(b)\setminus K$ we set
\[
g(x)\defeq \int_0^{\ima {\ts b}}|f^\prime(x+iy)|dy
\]
and if $x\in  K$ then we set $g(x)=0$. From the definition of $K$ we see that
\[
g(x)\leq 2^\ell |f^\prime(b)|\ima b\int_0^{\ima
  {\ts b}}\frac{|f^\prime(x+iy)|}{|f(x+iy) -f(b)|}dy  \] 
whenever $ x \in I ( b )$.  Using this inequality and Integrating   over 
$ I ( b ) $ we find that
\[
\int_{I(b)}g(x)dx\leq C2^\ell |f^\prime(b)|(\ima b)^2\leq
C^2  2^\ell d(f(b),\ar\Om)\ima b.
\]
So from weak-type estimates, if
\[
K^\prime\defeq\{x\in I(b): g(x)> 2^{2\ell}d(f(b),\ar \Om)\},
\]
then  

\eqn{4.25} \beq 
H^1(K^\prime)\leq C^2 2^{-\ell}\ima b,
 \eeq 
for some absolute constant $C$.
Using 
 (4.24) and (4.25)  we can fix $\ell$ to be a
large absolute number so that
\[  
H^1(K\cup K^{\prime}) <   10^{-100 }\ima b.
\]  
With $\ell$ thus fixed we put
 \[  E(b)\defeq I(b)\setminus(K\cup K^{\prime}) \]  and conclude that 
   Lemma  4.21 is valid.    
   $ \Box$ 
\\ \\

\noindent{\bf Lemma  4.26.} 
{\em Let   $ b, C^{\star}$ be as  in Lemma  4.21 and
put   $ c^* =    4 (C^{\star})^{2}$.  Given
$0<\delta<10^{-1000}$, let
$\delta_\star=e^{- c^* /\delta}.$ 
 Then, whenever
  $x\in E(b)$ there
is an interval $J = J ( x ) $ centered at $x$ with
 \eqn{4.27} \begin{equation}\label{eq:h1j}
2\delta_{\star}\ima b \leq H^{1} (J)\leq C\delta^{1/2} \ima b 
 \leq   \frac{\ima  b}{10000}  \end{equation}
(for some absolute constant $C$)
and a subset  $F = F ( x ) \subset J$
with $H^{1} (F) \geq  (1-10^{-100})H^1 (J)$ so that 
 \eqn{4.28} \begin{equation}\label{eq:f}
\int_0^{{\ts \delta_\star\ima b }}|f^\prime(t+iy)|dy\leq
\delta d(f(b),\ar\Om) \mbox{ \,  for every $ t \in F.$} 
\end{equation} } 
\\  

\noindent{\bf Proof of  Lemma 4.26:} 
Given $x\in E(b)$ put $ b' = x + i \ima b $ and  let  
$ y_1, 0<y_1<\ima b, $  be such that
\[
d(f(x+iy),\ar\Om)>\frac{\delta}{C^{\star}}d(f(b),\ar \Om)
\]
for $y_1<y<\ima b$, but
\[
d(f(\hat{b}),\ar\Om)=\frac{\delta}{C^{\star}}d(f(b),\ar \Om)
\]
where $\hat{b}\defeq x+iy_{1}$.  
By   (4.4),  Lemma  4.21, and conformal invariance of
hyperbolic distance,
\[
\log \frac{\ima b}{y_{1}}\leq \rho_{\Hh} (\hat{b},b')\leq 4Q_{\Om} (f
( \hat b ),f ( b' ))\leq \frac{4C^{\star}}{\delta d (f
(b),\ar \Om)}\int_{y_{1}}^{\ima {\ts b} }|f^{\prime} (x+iy)|dy\leq 
\frac{4 (C^{\star})^{2}}{\delta}, 
\]
i.e., $y_{1}\geq \delta_{\star}\ima b.$ 
Let $J=I (\hat{b})$ and $F=E (\hat{b})$.
Then by  Lemma  4.21, $H^1 (F)\geq (1-10^{-100})H^{1} (J)$
and for $t\in E (\hat{b})$
\[
\int_0^{{\ts \delta_\star\ima b}}|f^\prime(t+iy)|dy\leq C^{\star}d (f
(\hat{b}),\ar \Om)=\delta d (f (b),\ar \Om). 
\]
Notice also that,
\[
H^{1} (J)=2\ima \hat{b}\geq 2\delta_{\star}\ima b.
\]
On the other hand, elementary distortion theorems for univalent
functions   (see for example [GM05, ch 1, section 4]) and the fact
that $\hat{b}\in Q (b)$ yield   for
some absolute constant $ C_+  \geq 1 $ that  
 \[        \de/C^*  =  \frac{ d ( f ( \hat b ),  \ar \Om)}{ d  (
f ( b ) , \ar \Om ) }   \geq   \left ( \frac{\ima \hat b }{ C_+ \ima b }
\right)^2 .
\]  Thus   
(4.27), (4.28) are valid and the proof of  Lemma  4.26 is complete. 
 $ \Box $ \\ 

\noindent {\bf Lemma 4.29} { \em Let $ b, x \in E ( b ), J ( x ), F ( x
),  $ be as in Lemma  4.26 and set    $ \hat F = \bigcup_{x
\in  E ( b ) }  \, F ( x ). $  If $ L \subset I ( b) $ is an
interval with $ H^1 ( L ) \geq   {\ds \frac{ \ima b  }{ 100} }, $   then 
 \eqn{4.30}  \beq  H^1 (  E ( b ) \cap \hat F  \cap L ) \geq  
\frac{ \ima b }{1000} \, \, .
\eeq    Moreover, if   $    \{ \tau_1, \tau_2, \dots, \tau_m 
\} $ is a set of points in  $ I ( b ), $ then there exists  $
\tau_{m+1} $ in $  E ( b ) \cap \hat F  \cap L $ 
  with  
 \eqn{4.31} \beq      | f ( \tau_{m+1} ) - f ( \tau_j ) | \geq
  \frac{ d ( f ( b ), \ar \Om ) }{ 10^{10} \, \,  m^2 }  \mbox{ whenever }  
1 \leq j \leq m.
\eeq  }  \\
\noindent {\bf Proof of Lemma 4.29:}    Given an interval $ I $
let $ \la I $ be the interval with the same center as $ I $ and 
  $  \la  $ times its length.  Using Vitali, we see there exists   
$  \{  \hat  x_j  \} \subset     
  E ( b )  \cap \frac{1}{2} L  $  and  $ \{ J ( \hat x_j ) \} $
as in Lemma 4.26 such that 
  \[  E ( b )  \cap { \ts \frac{1}{2} } L   \subset
\bigcup_j \,  4  J ( \hat x_j ) \mbox{  and  the intervals $  \{ J ( \hat x_j )
\}   $  are pairwise disjoint.} \]

 Observe from  (4.27) that $  J ( \hat x_j ) \subset L $ for each
$ j. $   From this fact and  (4.27) we get  
\eqn{4.32}  \begin{eqnarray}   H^1 ( \hat F \cap L ) \geq 
 { \ds \sum_{j}  H^1 ( \hat F \cap J ( \hat x_j ) )  
   }  & \geq & ( 1 - 10^{ - 100} ) { \ds  \sum_j  H^1 ( J ( \hat x_j ) )     
 \geq     { \ts \frac{   1 - 10^{ - 100}  }{4} }     \, \sum_j   
 H^1 (  4 J ( \hat x_j ) )  } \notag \\ 
 &   \geq &   \frac{    1 - 10^{ - 100}  }{4}  \,  H^1 (  
 E ( b ) \cap  \frac{1}{2} L )  \geq   \ima b/ 900. 
 \end{eqnarray}    
 From (4.32) and (4.23) we conclude that  (4.30) is valid. 
To prove (4.31)  observe from  (4.30) and the Poisson integral
formula for $  \Hh $   that 
\eqn{4.33}  \beq  \om  (  E ( b ) \cap \hat F \cap L, b  )  \geq
10^{-4}  \eeq where $ \om ( \cdot, b )  $ denotes  harmonic measure on 
$ \Hh $ relative to $b.$   
Let  
\[  r = \sup_{ x \in E ( b ) \cap \hat F \cap  L } 
  \min  \{ | f ( x ) - f (  \tau_ j ) |, 1 \leq j \leq m \}.  \] 
Then  
\[ f ( E ( b ) \cap \hat F  \cap  L  ) \subset \bigcup_{j = 1}^m 
  \bar B ( f ( \tau_j ), r ). \] Using this fact,  (4.33), and invariance of  
harmonic measure under $ f, $ it follows that 
\eqn{4.34} \beq   10^{ - 4}  \leq  \sum_{j=1}^m   
 \ti \om ( \bar B ( f ( \tau_j ), r ), f ( b ) )  \eeq 
where $ \ti \om ( \cdot,  f ( b ) ) $ denotes harmonic measure in 
$ \Om $ relative to $ f ( b ). $  Finally we note  from the
Beurling projection theorem (see [GM05, ch 3,  Corollary 9.3])  that 
for each $ j, $ 
\[  \ti  \om ( \bar B ( f ( \tau_j ), r ), f ( b ) )  \leq  
     ( 4 / \pi )     \left(  \frac{r}{ d ( f ( b ), \ar \Om ) } 
\right)^{1/2}. \]   
    Using this inequality in (4.34) we conclude
that  (4.31) is true. The proof of Lemma 4.29 is now complete. 
$ \Box $ 

\section{ Proof of Theorem 2.} 
\subsection{Proof of (\ref{eq:silength}) and
(\ref{eq:w0})}\label{ssec:pfsilenghtw0}
Using Lemma 4.29 with $ b = a = is, $  
we deduce for given $ \de, 0 < \de  <  10^{ -
1000}, $  the existence of $ x_1, x_2, x_3  \in E ( a) $ 
with $  -  s  < x_1 < -  s/2,  
- \frac{1}{8}
 s  < x_3   < \frac{1}{8} s, $  and  $  \frac{1}{2} s  <
x_2  < s,  $ such that  
\eqn{5.1} \beq \int_{0}^{ \de_*  s}   | f' ( x_j + i  y ) | \,  dy  
\leq  \, \de d ( f ( a ), \ar \Om )    \mbox{ for } 1 \leq j
\leq 3,  \eeq   
 
\eqn{5.2} \beq   \min \{   | f ( x_1 ) - f ( x_3 ) |, | f ( x_2 ) - f
( x_3 ) | \} \geq  10^{ - 11}  d ( f ( a ), \ar \Om ).  \eeq 
\ 
As  earlier we let  $ \ti Q (a)$  be the shifted box whose
boundary in $ \Hh,    \xi, $ consists of the horizontal line segment from 
$ x_1 + is $ to $ x_2 + i s, $  and the vertical line segments
from $ x_j $ to $ x_j + i s, $ for $ j = 1, 2. $ Also  we
put $ \si = f ( \xi ) $ and note  from  $ x_j \in  E ( a ), j =
1, 2,  $  that $ (4.14)   $ is valid. Moreover, we let $w_i=f(x_i)$ for $i\in\{1,2,3\}$.
To construct $ \tau $ as defined  after (4.15), we  put $ t_1 =  x_3 $
and continue as outlined above (4.17). In general if  
 $ a_{k-1} = t_{k-1} + i s_{k-1}, $ we choose  $ t_k \in 
 E ( a_{k-1} )  $  so  that (4.17) holds with 
 $ s_k = \de_* s_{k-1}. $  This choice  is possible thanks to Lemma 4.29. 
 With  $ \la $ now defined note from
the argument following (4.17) that  $ x_0 = \lim_{t \rar 1} \la
(t) $ exists, $ | x_0 | < 1/4, $  and   that  
$  \tau  = f ( \la ) $ satisfies the
cigar condition in (4.16) for $  t \in [0,1). $   If  
$ w_0 = \tau (1),  $ then  
 using   (4.17), (4.18),   we see that   
\[   | w_3  - w_0 | \leq   \hat C  \de   d ( z_1, \ar
\Om )  \]  for some absolute constant $  \hat C. $  From this
inequality and  (5.2), it follows that  if      
 \eqn{5.3}   \beq    \de  =   \min ( 10^{-12}   \hat C^{ - 1},
10^{- 1000} ),    \eeq   
 then   \eqn{ 5.4} \beq   
\min \{ | w_0 - w_j |, j = 1, 2 \} \geq 10^{ - 12} d (  z_1, \ar
\Om ).  \eeq  
 With $ \de $ now fixed,  we see from (5.1) that  the part of $
\si, $ say $ \si_1, $ corresponding to   the vertical line
segments from $ x_j 
$ to $ x_j + i \de_* s,  \,  j = 1, 2, $  satisfies 
\eqn{5.5} \beq  d ( \si_1, w_0 )  \geq 10^{ - 13} d (  z_1,
\ar \Om ).  \eeq 
Using (4.4)  we also get  
\eqn{5.6}  \beq d ( \si\sem\si_1, \ar \Om ) \geq C^{-1} d ( z_1,
\ar \Om )   \eeq for some absolute constant $ C.  $   Combining
(5.5), (5.6), we obtain (4.15).    

\subsection{Proof of (\ref{eq:maxpr2})}\label{ssec:maxpr2}
\begin{figure}
%%%%%Figure 2
%\hspace{70pt}
\psfrag{a}{\large $a$}
\psfrag{Qa}{\large $\tilde{Q} (a)$}
\psfrag{l}{\large $\lambda$}
\psfrag{l1}{\large $\lambda_{1,1}$}
\psfrag{x1}{\large $x_{1}$}
\psfrag{xi}{\large $\xi$}
\psfrag{xi1}{\large $\xi_{1,1}$}
\psfrag{x0}{\large $x_{0}$}
\psfrag{x2}{\large $x_{2}$}
\psfrag{H}{\large $\Hh $}
\psfrag{b1}{\large $b_{1,1}$}
\psfrag{b2}{\large $b_{2,1}$}
\includegraphics{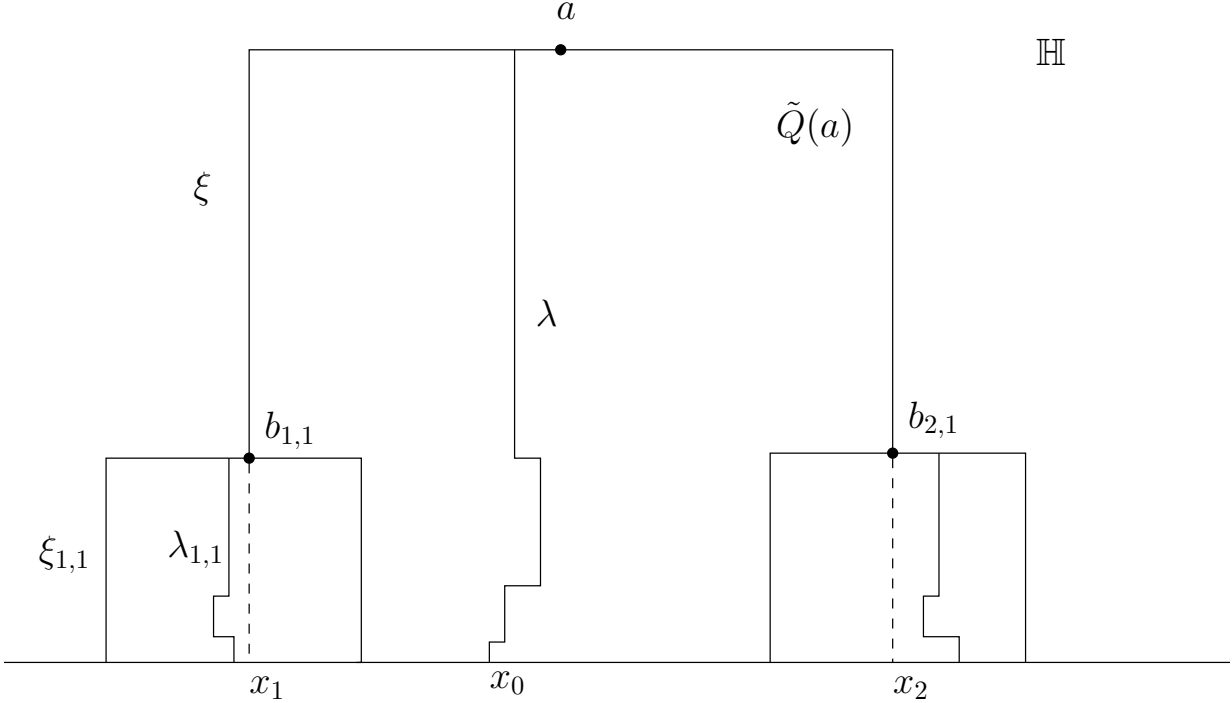}%
%\vspace{-.1in}

\caption{Domar-type recursion construction}\label{fig:domar}
\end{figure}

The proof of (4.13) is by contradiction. 
 Suppose $ u > A u ( z_1 ) $ on $ \si. $ We shall  obtain a 
contradiction   if 
$ A = A ( p )  $ is suitably large.  Our argument is based 
on  a recurrence  type scheme often attributed to Carleson -
Domar, see [C62], [D57], in the complex world, and to Caffarelli et. al., see [CFMS81], in the PDE world (see also [AS05] for references).    
 Given the shifted box
$\tilde{Q} (a)$ we  let   $b_{j,1} = x_j + i \de_* \ima a, j = 1,
2,  $  and
note that $ b_{j,1}, j = 1,2, $ are points on 
 the
vertical sides of 
$\tilde{Q} (a).$  These points  will  spawn two new boxes $\tilde{Q}
(b_{j,1})$, $j=1,2$, which in turn will each spawn two more new
boxes,  and so on. Without loss of generality, we focus
on  $\tilde{Q}
(b_{1,1}).$  This box is constructed in the same way as $\tilde{Q} (a)$  
and we also construct, using Lemma 4.29 once again,   
 a polygonal path $\lambda_{1,1}$ from $b_{1,1}$ to some point $x_{1,1}\in I
(b_{1,1}),$ so that $ \la_{1,1} $  is defined relative to $ b_{1,1} $ 
in the same way that $ \la $ was defined relative to $ a. $ 
 There is only one caveat. Namely, the path $\lambda_{1,1}$
is required to be contained in the half-plane $\{\rea z<\rea b_{1,1}
\}$, i.e., to stay entirely to the left of $b_{1,1}.$  
This extra caveat is easily achieved in view
of  Lemma 4.29. 
$ \la_{2,1} $  with endpoints, $ b_{2,1}, x_{2,1} $ is  constructed similarly,
 to lie in $ \{ \rea z > \rea b_{2,1} \}
$   (see Picture \ref{fig:domar}).

%%%%%%%%%%%%%%%%%%%%%%%%%%%%%%%%%%%%%%%%%%%%%%%%%%%%%%%%%%%%%%%%%%%%%%%%%%%%%%%%%%%%%% \begin{figure}
%%%%%Figure 2
%\hspace{70pt}

 %\psfrag{a}{\large $a$}
%
%
 %\psfrag{Qa}{\large $\tilde{Q} (a)$}
 %\psfrag{l}{\large $\lambda$}
  %\psfrag{l1}{\large $\lambda_{1,1}$}
   %\psfrag{x1}{\large $x_{1}$}
   %\psfrag{xi}{\large $\xi$}
  %\psfrag{xi1}{\large $\xi_{1,1}$}
  %\psfrag{x0}{\large $x_{0}$}
  %\psfrag{x2}{\large $x_{2}$}
  %\psfrag{H}{\large $\Hh $}
  %\psfrag{b1}{\large $b_{1,1}$}
  %\psfrag{b2}{\large $b_{2,1}$}
  %\includegraphics{domar.eps}%
%\vspace{-.1in}
%%%%%%%%%%%%%%%%%%%%%%%%%%%%%%%%%%%%%%%%%%%%%%%%%%%%%%%%%%%%%%%%%%%%%%%%%%%%%%%%%%%%%%%%%%%%%%%%%%%%% \caption{Domar-type recursion construction}\label{fig:domar} \end{figure}

%%%%%%%%%%%%%%%%%%%%%%%%%%%%%%%%%%%%%%%%%%%%%%%%%%%%%%%%%%%%%%%%%%%%%%%%%%%%%%%%%%%%%%%%%%%%%%%%%%%%%

Next, using the Harnack inequality we see that there exists $\Lambda$ such that
\eqn{5.7} \beq 
  u (f (z))\leq \Lambda u (  f (a)  )\mbox{ whenever }z=x+iy\in\xi,\ y\geq \delta_\ast\ima a.
\eeq 
 In particular, from  Harnack's inequality for $ u $ and
the fact that $ \de $ is now fixed in (5.3),  it is 
clear that  $ \La $ in (5.7)  can be chosen 
to depend only on $p$,  and hence  can also be used  in further
iterations.

By (5.7), the fact that $A>\Lambda$ and the maximum principle, we see that there exists a point 
$z\in \lambda_{1,1}\cup \lambda_{2,1}$  such that $U
(z)>AU (a)$. This is the reason why the paths $\lambda_{j,1}$ are
constructed outside the original box $\tilde{Q} (a)$.  First
suppose  $z\in
\lambda_{1,1}$. The
larger the constant $A$, the closer $z$ will be to 
$\re$. More precisely,
if $A>\Lambda^k$ then $\ima z\leq \delta_{\star}^{k}\ima a$, as
we see from (5.7) and inequalities analogous to (4.17)-(4.19).
 Arguing as in the display below (4.19), we find that   
\[
|f (z)-f (x_{1,1})|\leq C\delta^{k-1}d(f (b_{1,1}),\ar \Om).
\]
The argument now is similar to the argument showing the existence
of $ z^* $ at the end of subsection 4.1. 
Let $\xi_{1,1}$ be the boundary of $\tilde{Q} (b_{1,1})$ which is in $\Hh$
and let $\sigma_{1,1}=f (\xi_{1,1})$. Set $\rho_{1,1}\defeq d
(w_{0,1},\sigma_{1,1})$, where $w_{0,1}=f (x_{1,1})$. Then 
\[
B(w_{0,1},\rho_{1,1})\cap \Om \subset f (\tilde{Q} (b_{1,1})).
\]
So, by Lemma  2.2,
\[
u (f (z))\leq C\delta^{\alpha k}
\max_{\tilde{Q} (b_{1,1})}u\circ f.
 \,  \]
Choose $k$, depending only on $p$,  to be the least positive integer
such  that 
\[
C\delta^{\alpha k} <  \Lambda^{-1}.
\]
This choice of $k$ determines $A$ (say $A=2\Lambda^k$) which therefore
also depends only on $p$ (since $ \de $ is fixed in (5.3)).
With this choice of $A$ we have
\eqn{5.8} \beq 
\max_{\xi_{1,1}}U>\Lambda U (z)  > \Lambda A U (a).
\eeq 
 Since  $ U  ( b_{1,1} )  \leq  \La U (a) $  we see from  (5.8) 
 that we can now repeat the  above  argument with 
 $\tilde{Q} (b_{1,1})$ playing the role of $ \ti Q ( a ). $  That
is,  we find
$b_{1,2}$  on the vertical sides of $ \ti Q ( b_{1,1}) $ 
 with $\ima b_{1,2}=\delta_{\star}^{2}\ima a$  and  a box
$\tilde{Q} (b_{1,2})$ with boundary $\xi_{1,2}$ such that
\[
\max_{\xi_{1,2}}U > \Lambda^{2} A U (a)   \geq   A  U ( b_{1,2} ) .
\]
Continuing by induction we get a contradiction because $U=0$
continuously on $\re$.    If  $ z \in \la_{2,1}, $ we get a
contradiction by the same argument.  
 Thus,   there exists $A = A ( p ) \geq 1$  
 for which  (\ref{eq:maxpr2}) holds. The proof of Theorem 2 is
now complete. $ \Box $

\end{document}